\documentclass[11pt]{amsart}
\usepackage{amssymb}
\usepackage{amsmath}
\usepackage{mathabx}
\usepackage{amsthm}
\usepackage{amsmath,stmaryrd}
\usepackage{tikz-cd}
\usepackage{mathabx}
\usepackage{mathrsfs}
\usepackage{graphicx}
\usepackage{todonotes}
\usepackage{hyperref}
\usepackage[capitalize]{cleveref}
\usepackage{marginnote}
\usepackage{thmtools}
\usepackage[backend=biber, url=false, maxcitenames=10, maxbibnames=10]{biblatex}
\providecommand{\tightlist}{%
  \setlength{\itemsep}{0 pt}\setlength{\parskip}{0 pt}}

\DeclareMathOperator*{\bigdoublevee}{\bigvee\mkern-15mu\bigvee}

\newcommand{\edge}{\mathbin{E}}
\declaretheorem[numberwithin=section,name=Theorem]{thm}

\declaretheorem[name=Definition, style=definition, sibling=thm]{defn}
\declaretheorem[sibling=thm,name=Lemma]{lem}
\declaretheorem[sibling=thm]{proposition}
\declaretheorem[sibling=thm]{question}

\declaretheorem[sibling=thm]{example}
\declaretheorem[sibling=thm, name=Corollary]{cor}

\newcommand{\A}{{\mathcal{A}}}
\newcommand{\B}{{\mathcal{B}}}
\newcommand{\N}{\omega}
\newcommand{\mbf}[1]{\mathbf{#1}}
\newcommand{\mf}[1]{\mathfrak{#1}}

\renewcommand{\bar}[1]{\overline{#1}}
\newcommand{\restupa}[1]{P(\A)\restriction #1}
\DeclareMathOperator{\redu}{\leq_{p1}}
\renewcommand{\phi}{\varphi}

\newcommand{\oneone}{1{-}1}
\newbox\gnBoxA

\newdimen\gnCornerHgt
\setbox\gnBoxA=\hbox{$\ulcorner$}
\global\gnCornerHgt=\ht\gnBoxA
\newdimen\gnArgHgt
\def\Godelnum #1{%
\setbox\gnBoxA=\hbox{$#1$}%
\gnArgHgt=\ht\gnBoxA%
\ifnum \gnArgHgt<\gnCornerHgt \gnArgHgt=0pt%
\else \advance \gnArgHgt by -\gnCornerHgt%
\fi \raise\gnArgHgt\hbox{$\ulcorner$} \box\gnBoxA %
\raise\gnArgHgt\hbox{$\urcorner$}}

\hypersetup{
    colorlinks=true,
    }

\usepackage[margin=1.5in]{geometry}

\title{Relations enumerable from positive information}

\author{Barbara F. Csima}
\address{Department of Pure Mathematics, University of Waterloo}
\email{csima@uwaterloo.ca}
\author{Luke MacLean}
\address{Department of Mathematics and Statistics, University of New Brunswick}
\email{luke.maclean@unb.ca}

\author{Dino Rossegger}
\address{Institute of Discrete Mathematics and Geometry, Technische Universit\"at Wien}
\email{dino.rossegger@tuwien.ac.at}

\thanks{The first author is supported by an NSERC Discovery grant. The work of the third author was supported by the European Union's Horizon 2020 Research and Innovation Programme under the Marie Sk\l{}odowska-Curie grant agreement No. 101026834 — ACOSE. This project received support from the Austrian Agency for International Cooperation in Education and Research under grant WTZ-BG11/2019.}

\begin{document}
\begin{abstract}
    We study countable structures from the viewpoint of enumeration reducibility. Since enumeration reducibility is based on only positive information, in this setting it is natural to consider structures given by their positive atomic diagram -- the computable join of all relations of the structure. Fixing a structure $\A$, a natural class of relations in this setting are the relations $R$ such that $R^{\hat \A}$ is enumeration reducible to the positive atomic diagram of $\hat \A$ for every $\hat \A\cong \A$ -- the \emph{relatively intrinsically positively enumerable (r.i.p.e.)} relations. We show that the r.i.p.e.\ relations are exactly the relations that are definable by $\Sigma^p_1$ formulas, a subclass of the infinitary $\Sigma^0_1$ formulas. We then introduce a new natural notion of the jump of a structure and study its interaction with other notions of jumps.
\end{abstract}
\maketitle
While in computable structure theory countable structures are classically studied up to Turing reducibility, researchers have successfully used enumeration reducibility to both contribute to the classical study and develop a beautiful theory on its own. One example of a contribution to classical theory is Soskov's work on degree spectra and co-spectra~\cite{soskov2004} which allowed him to show that there is a degree spectrum of a structure such that the set of degrees whose $\omega$-jump is in this spectrum is not the spectrum of a structure~\cite{soskov2013a}. Another example is Kalimullin's  study of reducibilities between classes of structures~\cite{kalimullin2012} where he studied enumeration reducibility versions of Muchnik and Medvedev reducibility between classes of structures. This topic has also been studied by Knight, Miller, and Vanden Boom~\cite{knight2007} and Calvert, Cummins, Knight, and Miller~\cite{calvert2004a}. There appears to be a rich theory on these notions with interesting questions on the relationship between enumeration reducibility and the classical versions. In this article we develop a novel approach to study relations that are enumeration reducible to every copy of a given structure.

As mentioned in the abstract, since enumeration reducibility is based only on positive information, it is natural to consider structures given by their positive atomic diagram in this setting. This is quite reasonable. When doing computability theory and looking at an infinite object, we often view it as being revealed stage by stage. It is natural to find out only about positive information, and our approach is general enough to handle negative information by considering the associated structure in the language including relation symbols for the cosets of the relations in the original language.

More formally, say that given a countable relational structure $\A$, a relation $R$ is \emph{relatively intrinsically positively enumerable} (\emph{r.i.p.e.}) in $\A$ if for every copy $\B$ of $\A$ and every enumeration of the basic relations on $\B$, we can compute an enumeration of $R^\B$.
We obtain a syntactic classification of the r.i.p.e.\ relations using the infinitary logic $L_{\omega_1\omega}$. \cref{ripe} shows that the r.i.p.e.\ relations on a structure are precisely those that are definable by computable infinitary $\Sigma_1^0$ formulas in which neither negations nor implications occur.

This classification is motivated by a related classification of the r.i.c.e.\ relations. A relation $R$ on a structure $\A$ is relatively intrinsically computably enumerable (r.i.c.e.) if for every $\B\cong \A$, $R^\B$ is c.e.\ relative to the atomic diagram of $\B$.
The following classification of r.i.c.e.\ relations is a special case of a theorem by Ash, Knight, Manasse and Slaman~\cite{ash1989} and, independently, Chisholm~\cite{chisholm1990}: A relation is r.i.c.e.\ in a structure $\A$ if and only if it is definable in $\A$ by a computable $\Sigma^0_1$ formula in the logic $L_{\omega_1\omega}$.

Since this result there has been much work on r.i.c.e.\ relations and related concepts. For a summary, see Fokina, Harizanov, and Melnikov~\cite{fokina2014}. One particularly interesting generalization of r.i.c.e.\ relations is due to Montalb\'an~\cite{montalban2012}. He extended the definition of r.i.c.e.\ relations from relations on $\omega^n$ to $\omega^{<\omega}$ and to sequences of relations, obtaining a classification similar to the one given in~\cite{ash1989,chisholm1990}. This extension allows the development of a rich theory such as an intuitive definition of the jump of a structure. For a complete development of this theory see Montalb\'an~\cite{montalban2021a}. The main goal of this article is to develop a similar theory for r.i.p.e.\ relations.
\begin{defn}\label{def:posdiagram}
Let $\A$ be a countable structure in relational language $L$ with universe $\omega$. We let $P(\A)$ be the set $=^\A\oplus \neq^\A\oplus\bigoplus_{R_i\in L}R_i^\A$, which we call the \emph{positive diagram of $\A$}.
\end{defn}
The mindful reader might notice that we include inequality in the positive atomic diagram. The reason for this is that for enumerations of $\A$ that are not injective $=^\A$ will become an equivalence relation instead of equality and we want to take the quotient under this equivalence relation to be able to enumerate an isomorphic copy of $\A$ from this enumeration. Having $\neq$ in the diagram allows us to do this effectively. The positive diagram is Turing equivalent to the standard definition of the atomic diagram of a structure, which can be viewed as the set $=^\A\oplus \neq^\A\oplus\bigoplus_{R_i\in L}R_i^\A\oplus \bar{R_i}^\A$. Now, a relation $R\subseteq\omega^{<\omega}$ is r.i.c.e.\ if for every copy $\B\cong \A$, $R^\B$ is c.e.\ in $D(\B)$. 
The r.i.p.e.\ relations are the natural analogue to the r.i.c.e.\ relations for enumeration reducibility. Recall that a set of natural numbers $A$ is enumeration reducible to $B$, $A\leq_e B$ if there exists a c.e.\ set $\Psi$ consisting of pairs $\langle x,D\rangle$ where $D$ is a finite set under some fixed coding such that
\[ A=\Psi^B=\{x: D\subseteq B \land \langle x,D \rangle \in \Psi\}. \]
Enumeration reducibility allows us to formally define the notion of a r.i.p.e.\ relation.
\begin{defn}\label{ripedefn}
    Let $\A$ be a structure. A relation $R \subseteq A^{< \omega}$ is \emph{relatively intrinsically positively enumerable} in $\A$, short r.i.p.e., if, for every copy $(\B,R^\B)$, of $(\A, R)$, $R^\B\leq_e P(\B)$. The relation is \textit{uniformly relatively intrinsically positively enumerable} in $\A$, if the above enumeration reducibility is uniform in the copies of $\A$, that is, if there is a fixed enumeration operator $\Psi$ such that $R^\B = \Psi^{P(\B)}$ for every copy $\B$ of $\A$.
\end{defn}
The study of the computability theoretic properties of structures with respect to enumeration reducibility is an active research topic; see Soskova and Soskova~\cite{soskova2017} for a summary of results in this area. 
\subsection{Structure of the paper}
In \cref{sec:ripe} we show that a relation is r.i.p.e.\ in a structure $\A$ if and only if it is definable by a $\Sigma_1^p$ formula in the language of $\A$, that is, a computable infinitary $\Sigma_1^0$ formula in which neither negations nor implications occur. We continue by defining a notion of reducibility between r.i.p.e.\ relations and exhibiting a complete relation. Towards the end of \cref{sec:ripe} we study sets of natural numbers r.i.p.e.\ in a given structure and show that these are precisely those sets whose degrees lie in the structures co-spectrum, the set of enumeration degrees below every copy of the structure.

\cref{sec:jumps} is devoted to the study of the positive jump of a structure. We define the positive jump and study its degree theoretic properties. The main result of this section is that the enumeration degree spectrum of the positive jump of a structure is the set of jumps of enumeration degrees of the structure. The main tool to prove this result are r.i.p.e.\ generic enumerations which are studied at the beginning of the section. To finish the section we compare the enumeration degree spectra obtained by applying various operations on structures such as the positive jump, the classical jump and the totalization.

%


\section{First results on r.i.p.e.\ relations}\label{sec:ripe}

In our proofs we will often build structures in stages by copying increasing pieces of finite substructures of a given structure $\A$. The following definitions will be useful for this.
\begin{defn}
   Given an $\mathcal L$-structure $\A$ and $\overline a\in A^{<\omega}$, let $\restupa{\overline a}$ denote the positive diagram of the substructure of $\A$ with universe $\overline a$ in the restriction of $\mathcal L$ to the first $|\overline a|$ relation symbols.
\end{defn}

\begin{defn}\label{partial}
    Let $\A$ be an $\mathcal L$-structure and $\overline{a}= \langle a_0, \dots , a_{s-1} \rangle\in A^s$. The set $P_{\A}(\overline a)$ is the pullback of $\restupa{\overline a}$ along the index function of $\overline a$, i.e., 
    \[ \langle i,n_0,\dots, n_{r_i}\rangle\in P_{\A}(\overline a) \Leftrightarrow \langle i, a_{n_0},\dots a_{n_{r_i}}\rangle \in \restupa{\overline a}.\]
\end{defn}
The main feature of \cref{partial} is that if we approximate the positive diagram of a structure $\A$ in stages by considering larger and larger tuples, i.e., $\lim_{s\in\omega} \restupa{\overline a_s}=P(\A)$, then the limit of $P_{\A}(\overline a_s)$ gives a structure isomorphic to $\A$. This fact is useful in constructions.  

We denote by $\Psi_e$ the $e^{th}$ enumeration operator in a fixed computable enumeration of all enumeration operators and by $\Psi_{e,s}$ its stage $s$ approximation. Without loss of generality we make the common assumption that $\Psi_{e,s}$ is finite and does not contain pairs $\langle n,D\rangle>s$. Notice that $\Psi_{e,s}$ itself is an enumeration operator.
In our proofs we also frequently argue that a set $A$ is enumeration reducible to a set $B$ by using a characterization of enumeration reducibility due to Selman~\cite{selman1971a}. 
\begin{thm}[Selman~\cite{selman1971a}]
    For any $A, B \subseteq \omega$
    \begin{equation*}
        A \leq_e B \quad \text{ iff } \quad \forall X[B \text{ is c.e. in } X \Rightarrow A \text{ is c.e. in } X]
    \end{equation*}
\end{thm}
We refer the reader to Cooper~\cite{cooper2003} for a proof of this result and further background on enumeration reducibility and enumeration degrees.

Notice that by \cref{ripedefn} $(\A,R)$ is technically not a first order structure as $R\subseteq \omega^{<\omega}$. We may however still think of it as a first order structure in the language expanded by relation symbols $(Q_i)_{i\in\omega}$, each $Q_i$ of arity $i$, where
$Q_i^\A=\{ \bar a \in A^{n}: \bar a\in R^\A\}$. The positive diagrams $P(\A,R^\A)$ and $P(\A,(Q_i^\A)_{i\in\omega})$ are enumeration equivalent.
\subsection{Examples of r.i.p.e.\ relations}
\begin{example}
Let $\mathcal R$ be a countable ring. Then the relation 
$P=\{ \bar a: \exists x [a_0+a_1x+\dots+a_kx^k=0]\}$
holding of the polynomials with a root in $\mathcal R$ is uniformly  r.i.p.e. \end{example}
\begin{proof}
Enumerate the graphs of $+$,$\cdot$, $=$ and $0$ of a  copy $\mathcal{\hat R}$ of $\mathcal R$ and whenever you see $\bar ab$ such that $a_0+a_1b+\dots+ a_kb^k=0$ enumerate $\bar a$. Clearly this gives an enumeration of $P^\mathcal{\hat R}$ from an enumeration of any copy of $\mathcal R$.
\end{proof}
\begin{example}\label{grapheg}
The reachability relation $reach(x,y)$ on graphs is uniformly r.i.p.e.
\end{example}
\begin{proposition}
There is a graph $G$ in the language $(=,\neq,E)$ such that $\neg E$ is not r.i.p.e.\ in $G$.
\end{proposition}
\begin{proof}
Let $K$ be the halting set and let $G$ be the following graph: 
\begin{enumerate}
\tightlist
    \item For every $n\in\omega$, $G$ contains a vertex $v_n$ and an $(n+1)$-cycle with least vertex $c_n$ and  an edge $v_n\edge c_n$,
    \item if $\langle x,y\rangle\in K$, then $v_{2x+1}\edge v_{2y+2}$.
\end{enumerate}
Then it is not hard to see that from any enumeration of $K$ we can enumerate $P(G)$. Notice that the set $C=\{ (v_{2x+1},v_{2y+2}): (x,y)\in \omega\times\omega\}$ is r.i.p.e. If $\neg E$ was r.i.p.e., then so would be $C\cap \neg E$, i.e., $C\cap \neg E\leq_e P(G)$. But clearly $C\cap\neg E \geq_e \bar K$ and since we have $K\geq_e P(G)$ this would imply $K\geq_e \bar K$, a contradiction. 
\end{proof}
\subsection{A syntactic characterization}
The main purpose of this section is to show that being r.i.p.e.\ in a structure $\A$ is equivalent to being definable by infinitary formulas in $\A$ of the following type.
\begin{defn}
    A \emph{positive computable infinitary $\Sigma_1^p$ formula} is a formula of the infinitary logic $L_{\omega_1\omega}$ of the form
    \begin{equation*}
        \varphi(\overline{x}) \equiv \bigdoublevee_{i \in I} \exists \overline{y}_i \varphi_i (\overline{x},\overline{y}_i)
    \end{equation*}
    where each $\varphi_i$ is a finite conjunct of atomic formulas, and the index set $I$ is computably enumerable.
\end{defn}
Notice that the above definition includes all c.e.\ disjunctions of finitary $\Sigma^0_1$ formulas without negation and implication symbols, as every such formula is equivalent to a finite disjunction of conjunctions with each existential quantifier occurring in front of the conjunctions.
\begin{defn}
    A relation $R\subseteq \N^{<\N}$ is \emph{$\Sigma^p_1$-definable with parameters $\bar c$} in a structure $\A$ if there exists a uniformly computable sequence of $\Sigma^p_1$ formulas $(\phi_i(x_1,\dots,x_i, y_1,\dots,y_{|c|}))_{i\in\omega}$ such that for all $\bar a\in \N^{<\N}$ 
    \[ \bar a\in R \Leftrightarrow \A\models  \phi_{|\bar a|}(\bar a, \bar c).\]
\end{defn}
Our goal is to show that a relation $R$ is r.i.p.e.\ in a structure $\A$ if and only if it is $\Sigma^p_1$-definable over some parameter. We will do this by using forcing to build a generic copy of $\A$ and then read the syntactic definition of this copy. The right notion of genericity for this purpose is the following.
\begin{defn}\label{def:ripegeneric}
 Let $A^* = \{\sigma \in A^{<\omega} : (\forall i \neq j < |\sigma|)[\sigma(i) \neq \sigma(j)]\}$. We say that $\gamma \in A^*$ \emph{decides} an upwards closed subset $R\subseteq A^*$ if $\gamma \in R$ or $\sigma \not \in R$ for all $\sigma \supseteq \gamma$. A $\oneone$ function $g:\omega \rightarrow \A$ is a \emph{r.i.p.e.-generic enumeration}, if for every r.i.p.e. subset $R \subseteq A^*$ there is an initial segment of $g$ that decides $R$. We say that $\B$ is a \emph{r.i.p.e.-generic presentation} of $\A$ if it is the pull-back along a r.i.p.e.-generic enumeration of $P(\A)$.
\end{defn}
The existence of r.i.p.e.\ generics follows from the Baire category theorem. Let us give a quick hands-on construction. Given a structure $\A$ and an enumeration $(R_e)_{e\in\omega}$ of all r.i.p.e.\ subsets of $A^*$ we build an sequence $\bar p_0=\emptyset\subset \bar p_1\dots$ with $\lim p_s=g$. At even stages $2s$ we guarantee that $g$ is onto by defining $\bar p_{2s}={\bar p_{2s-1}}^\smallfrown x$ where $x\in A$ is least such that $x$ is not in the range of $p_{2s-1}$. At odd stages $2s+1$ we check whether there is $\bar q\supseteq\bar p_{2s}$ with $\bar q\in R_s$. If so, let $\bar p_{2s+1}=\bar q$, otherwise $\bar p_{2s+1}=\bar p_{2s}$. It is easy to see that $g=\lim p_s$ decides every upwards closed r.i.p.e.\ subset of $A^*$ and is thus a r.i.p.e.\ generic enumeration. Furthermore, note that if $R_e$ is dense in $A^*$, then there is some $p\subset g$ forcing that $g$ meets $R_e$ and thus $\langle 0,\dots,|p|\rangle\in R_e^{g^{-1}(\A)}$. We will study properties of r.i.p.e.\ generics in \cref{sec:ripegenerics}. For now, let us obtain a syntactic characterization of r.i.p.e.\ relations.
\begin{thm}\label{ripe}
    Let $\A$ be a structure and $R \subseteq A^{< \omega}$ a relation on it. Then the following are equivalent:
    \begin{enumerate}\tightlist
        \item[(i)] $R$ is relatively intrinsically positively enumerable in $\A$,
        \item[(ii)] $R$ is $\Sigma^p_1$-definable in $\A$ with parameters.
    \end{enumerate}
\end{thm}
\begin{proof} Assuming (ii), there is a uniformly computable sequence $(\phi_i(\overline x,\overline z))_{i\in\omega}$ of $\Sigma_1^p$ formulas where each $\phi_i$ is of the form $\bigdoublevee_{j \in \omega} \exists \overline{y}_{i,j} \psi_{i,j} (\overline{x},\overline{y}_{i,j},\overline{z})$ with the property that for every $\B\cong \A$ there is a tuple $\bar c\in \omega^{|\bar z|}$ such that for all $i\in\omega$ and $\bar a\in\omega^{i}$, $\B\models \varphi_i(\bar a,\bar c)$ if and only if $\bar a \in R^{\B}$.
Recall that each $\psi_{i,j}$ is a finite conjunction of atomic formulas, i.e., $\psi_{i,j}=\theta_1(\overline x, \overline y_{i,j},\overline z)\land \dots \land \theta_n(\overline x, \overline y_{i,j},\overline z)$ for some $n\in \omega$. For $\theta(\overline x)$ an atomic formula, let $\ulcorner \theta(\overline a)\urcorner$ be the function mapping $\theta(\overline a)$ to its code in the positive diagram of a structure. For example, if $\theta(\overline x)=R_i(x_3,x_5)$, then $\ulcorner \theta(\overline a)\urcorner=\langle i+2, \langle a_3,a_5\rangle\rangle$ for $\overline a\in \omega^{<\omega}$. Consider the set $X^{\overline a, \overline b,\overline c}_{i,j}=\{\ulcorner \theta_k(\overline a,\overline b,\overline c)\urcorner: k<n\}$. Clearly, $X^{\overline a,\overline b,\overline c}_{i,j}\subseteq P(\A)$ if and only if $\A\models \psi_{i,j}(\overline a,\overline b,\overline c)$ for any $L$-structure $\A$. We define an enumeration operator $\Psi$ by enumerating all pairs 
$\langle \overline a,X_{i,j}^{\overline a,\overline b,\overline c}\rangle$ into $\Psi$. It follows that 
\begin{align*}
    \overline a\in \Psi^{P(\B)}&\Leftrightarrow \exists \langle \overline a,X_{|\overline a|,j}^{\overline a, \overline b, \overline c}\rangle \in \Psi \land X_{|\overline a|,j}^{\overline a, \overline b\overline c}\subseteq P(\B)\\
    &\Leftrightarrow \B\models \exists \overline y_{|\overline a|,j}\psi_{|\overline a|,j}(\overline a, \overline y_{|\overline a|,j},\overline c)\\
    &\Leftrightarrow \B\models \phi_{|\overline a|}(\bar a, \bar c)
\end{align*}
and thus $R$ is r.i.p.e.

Let $g$ be a r.i.p.e.\ generic enumeration of $\A$ produced as in the paragraph above the theorem and let $\B=g^{-1}(\A)$. Our goal is to produce a $\Sigma^p_1$ definition of $R^\B$, as $\B\cong \A$ we get that this is then also a definition for $R$. 
%
Towards this, consider the set
\begin{equation*}
    Q_e = \{ \overline{q} \in A^* : \exists l, j_1 ,\dots j_l < |\overline{q}| \, [\langle j_1,\dots, j_l \rangle \in \Psi_e^{P_\A(\overline{q})} \text{ and } \langle q_{j_1},\dots, q_{j_l} \rangle \not \in R] \}.
\end{equation*}
We have that $g$ meets $Q_e$ if and only if $R^\B\neq \Psi_e^{P_\B}$. By our assumption that $R$ is r.i.p.e., there is $e_0$ such that $\Psi^{P(\B)}_{e_0}=R^\B$, and thus $g$ does not meet $Q_{e_0}$.

We will use this to give a $\Sigma^p_1$ definition of $R$ with parameters $\bar p_{s}$, where $\bar p_{s}$ is the stage in the construction of $g$ at which we decide $Q_{e_0}$.

Notice that if there is some $\bar q \supseteq \bar p_s$ and sub-tuple $\langle q_{j_1} , \dots, q_{j_l} \rangle$ such that $\langle j_1,\dots , j_l \rangle \in \Psi_e^{P_\A(\overline{q})}$ then we must have $\langle q_{j_1} , \dots, q_{j_l} \rangle \in R$ or else we will have that $\bar q \in Q_{e_0}$. We now show that $R$ is equal to the set
\begin{align*}
    S = \{ \langle q_{j_1},\dots , q_{j_l} \rangle \in A^* : \text{ for some } \bar q \in A^{< \omega} \text{ and } l , j_1, \dots , j_l < |\bar q | \text{ satisfying } \bar q \supseteq \bar p_s \\ \text{ and } \langle j_1 ,\dots , j_l \rangle \in \Psi_e^{P_\A(\bar q)}\}
\end{align*}
By the previous paragraph, $S \subseteq R$. If $\bar a \in R$ then there are indices $j_1,\dots, j_{|\bar a|}$ such that $\bar a = \langle g(j_1), \dots g(j_{|\bar a|}) \rangle$ and so if we take a long enough initial segment of $g$ it will witness the fact that $\bar a \in S$. Fix an enumeration $(\phi_i^{at})_{i\in\omega}$ of all atomic formulas where without loss of generality the free variables of $\phi_i^{at}$ are a subset of $\{x_0,\dots, x_{i}\}$. The following is a $\Sigma^p_1$ definition of $S$
\[
    \bigdoublevee_{ C \subset_{\text{fin}} \omega} \bigvee_{\langle j_1,\dots , j_{|\bar a|} \rangle \in W_e^C} \exists\, \bar q \supseteq \bar p_s \, [\langle q_{j_1}, \dots ,q_{j_{|\bar a|}} \rangle = \bar a \, \wedge \, \bigwedge_{i \in C} [\varphi^{at}_i]\frac{q_0}{x_0}\cdots\frac{q_{|\bar q|}}{x_{|\bar q|}} ]
\]
where the latter half of the formula is simply saying that $C \subseteq P_\A(\bar q)$.
\end{proof}

\begin{cor}\label{uripe}
    Let $\A$ be a structure and $R \subseteq A^{< \omega}$ be a relation on it. Then the following are equivalent:
    \begin{enumerate}\tightlist
        \item[(i)] $R$ is uniformly relatively intrinsically positively enumerable in $\A$,
        \item[(ii)] $R$ is $\Sigma_1^p$-definable in $\A$ without parameters. 
    \end{enumerate}
\end{cor}
\begin{proof}
    For $(i) \Rightarrow (ii)$ we mimic the proof of \cref{ripe}. Let $\Psi_e$ be the fixed enumeration operator such that $R^{\B} = \Psi_e^{P(\B)}$ and $Q_e$ as above. Note that no $\bar q$ can be in $Q_e$ and so we mimic the construction of our set $S$ with $\bar p_s$ being the empty tuple.
    
    For $(ii) \Rightarrow (i)$ we again mimic the same direction in \cref{ripe} excluding the parametrizing tuple $\bar c$ to make the process uniform.
\end{proof}
\subsection{R.i.p.e.\ completeness}
Similar to the study of computably enumerable sets we want to investigate notions of completeness for r.i.p.e.\ relations on a given structure. Before we obtain a natural example of a r.i.p.e.\ complete relation we have to define a suitable notion of reduction.
\begin{defn}
Given a structure $\A$ and two relations $P,R\subseteq A^{<\omega}$, we say that $P$ is \emph{positively intrinsically one reducible} to $R$, and write $P\redu R$, if for all $\B\cong \A$ $P(\B)\oplus P^\B\leq_1 P(\B)\oplus R^\B$.
\end{defn}
\begin{proposition}
    Positive intrinsic one reducibility ($\redu$) is a reducibility.
\end{proposition}
\begin{proof}
    Let $\A$ be a structure with relations $P, Q, R \subseteq A^{<\omega}$. It is easy to see that $\redu$ is reflexive, since for any structure $\B \cong \A$ we have that $P(\B) \oplus P^\B \leq_1 P(\B) \oplus P^\B$, which means $P \redu P$. To see that it is transitive assume that $P\redu R$ and $R \redu Q$ and let $\B \cong \A$. By assumption $P(\B) \oplus P^\B \leq_1 P(\B) \oplus R^\B \leq_1 P(\B) \oplus Q^\B$ and so $P \redu Q$.
\end{proof}
\begin{defn}
    Fix a structure $\A$. A relation $R\subseteq A^{<\omega}$ is \emph{r.i.p.e.\ complete} if $R$ is r.i.p.e.\ and for every r.i.p.e.\ relation $P$ on $\A$, $P\redu R$.
\end{defn}
The most natural way to obtain a complete set is to follow the construction of the Kleene set in taking the computable join of all r.i.p.e.\ sets. The result is a relation $R\subseteq \omega\times A^{<\omega}$ which can be seen as a uniform sequence of r.i.p.e.\ relations in the sense that there is a enumeration operator $\Psi$ such that $({\Psi^{\A}})^{[i]}=R_{i}$. 
For any isomorphic copy of $\A$, the $n$th slice of its Kleene set should correspond to the $n$th r.i.p.e.\ relation. In order to accomplish this we use the following coding.
Given a relation $R \subseteq \omega \times A^{<\omega}$, we can identify it with a subset $R'\subseteq A^{<\omega}$ as follows. For any two elements $b,c \in A$ let
\[
\overbrace{b\dots b}^{i\times}c\bar a \in R'\iff\langle i, \bar a \rangle \in R.  
\]
One can now easily see that $R$ is a uniform sequence of r.i.p.e.\ relations if and only if the so obtained relation $R'$ is uniformly r.i.p.e.

We can now give a natural candidate for a r.i.p.e.\ complete relation.
\begin{defn}
    Let $\phi_{i,j}^{\Sigma^p_1}$ be the $i${th} formula with free variables $x_1,\dots, x_j$ in a computable enumeration of all $\Sigma^p_1$ formulas.
    The \emph{positive Kleene predicate} relative to $\A$ is $\vec K^\A_p = \left(K^\A_{i}\right)_{i\in\omega}$, where
    \begin{equation*}
      K^{\A}_{i} = \bigcup_{j\in\omega}\{ \bar a \in A^{j} : \A \models \varphi^{\Sigma_1^p}_{i,j}(\bar a)\}.
    \end{equation*}
\end{defn}
Notice that we defined the positive Kleene predicate as a sequence of relations instead of a single relation. It is slightly more convenient as we do not have to deal with coding, but we could write it as a single relation on $\omega \times A^{<\omega}$ and code it as above. Another alternative definition would be to break down the sequence even further and let the Kleene predicate be the sequence 
\[ (\{\overline a: \A\models \phi_{i,j}^{\Sigma_1^p}(\overline a)\})_{\langle i,j\rangle\in\omega}\]
so that $(\A,\vec K_p^\A)$ is a first order structure. However, as all of these definitions are computationally equivalent these distinctions are irrelevant for our purpose.


\begin{proposition} The positive Kleene predicate $\vec K_p^\A$ is uniformly r.i.p.e., and r.i.p.e.\ complete.
\end{proposition}
\begin{proof}
    Since $\vec K^\A_p$ is defined in a $\Sigma_1^p$ way without parameters we
    can use \cref{uripe} to see that it is uniformly relatively intrinsically positively enumerable. Let $R$ be any relation on $\A$ of arity $a_R$. Notice that $R^\B$ is trivially $\Sigma_1^p$ definable, and so there is a formula $\varphi_{i,a_R}^{\Sigma_1^P}$ such that $\bar a \in R^\B \Leftrightarrow \B \models \varphi_{i,a_R}^{\Sigma_1^p}(\bar a)$. This shows that $P(\B) \oplus R^\B \leq_1 P(\B) \oplus \vec K_p^\B$.
\end{proof}

\subsection{R.i.p.e.\ sets of natural numbers}
Our above discussion of sequences of r.i.p.e.\ relations allows us to code sets of natural numbers as r.i.p.e. relations. 
\begin{defn}
    A set $X\subseteq\omega$ is \emph{r.i.p.e.\ in a structure $\A$} if the relation $R_X=X\times \emptyset$ is a r.i.p.e.\ relation, i.e., if the following relation is r.i.p.e.:
    \[ \{\overbrace{b\dots b}^{i\times}c: i\in X, b,c\in A\}\]
\end{defn}
A natural question is which sets of natural numbers are r.i.p.e.\ in a given structure. One characterization can be derived directly from the definitions: The sets $X$ such that $X\leq_e P(\B)$ for all $\B\cong \A$. Another one can be given using co-spectra, a notion defined by Soskov~\cite{soskov2004}. Intuitively, the co-spectrum of a structure $\A$ is the maximal ideal in the enumeration degrees such that every member of it is below every copy of $\A$. More formally.
\begin{defn}\label{def:cospectra}
    The \emph{co-spectrum} of a structure $\A$ is
        \[ Co(\A)=\bigcap_{\B\cong \A}\{ \mathbf d : \mathbf d\leq deg_e(P(\B))\}.\]
\end{defn}
Let us point out that Soskov's definition of co-spectra appears to be different from ours. We will prove in \cref{sec:jumps} that the two definitions are equivalent. This definition also agrees with the definition of co-spectra on Turing degrees when we restrict our attention to total degrees.

\begin{defn}
    A set $A$ is said to be \emph{total} if $A \equiv_e A \oplus \bar A$. An enumeration degree is said to be total if it contains a total set, and a structure $\A$ is total if $P(f^{-1}(\A))$ is a total set for every enumeration $f$.
\end{defn}

Given a tuple $\bar a$ in some structure $\A$ let $\Sigma_1^p\text{-}tp_\A(\bar a)$ be the set of positive finitary $\Sigma_1^0$ formulas true of $\A$. The equivalence of \cref{it:ripe} and \cref{it:enum} in \cref{thm:ripesets} is the analogue to a well-known theorem of Knight~\cite{knight1986} for total structures.

\begin{thm}\label{thm:ripesets}
    The following are equivalent for every structure $\A$ and $X\subseteq \omega$.
    \begin{enumerate}\tightlist
        \item\label{it:ripe} $X$ is r.i.p.e.\ in $\A$
        \item\label{it:cosp} $deg_e(X)\in Co(\A)$
        \item\label{it:enum} $X$ is enumeration reducible to $\Sigma_1^p\text{-}tp_\A(\bar a)$ for some tuple $\bar a\in A^{<\omega}$.
    \end{enumerate}
\end{thm}
\begin{proof}
    If $deg_e(X)\in Co(\A)$, then for all $\B\cong \A$, $X\leq_e P(\B)$. Given an enumeration of $X$, enumerate $b^nc$ into $R_X$ for all elements $b,c\in B$ whenever you see $n$ enter $X$. The relation $R_X$ clearly witnesses that $X$ is r.i.p.e.\ in $\A$. This shows that \cref{it:cosp} implies \cref{it:ripe}. On the other hand if $X$ is r.i.p.e.\ in $\A$, then given any $B\cong \A$ and an enumeration of $R_X^\B$ build a set $S$ by enumerating $n$ into $S$ whenever you see $b^n c$ enumerated into $R_X^\B$ for any two elements $b,c\in B$. Clearly $n\in S$ if and only if $n\in X$ and thus \cref{it:ripe} implies \cref{it:cosp}.
    
    To see that \cref{it:enum} implies \cref{it:cosp} assume that $X$ is enumeration reducible to the positive $\Sigma_1$ type of a tuple $\bar b$ in any copy $\B$ of $\A$. As the $\Sigma_1^p\text{-}tp_\B(\bar b)$ is enumeration reducible to $P(\B)$, by transitivity $X\leq_e P(\B)$ for every $\B\cong \A$ and thus $deg_e(X)\in Co(\A)$.
    At last, we show that \cref{it:ripe} implies \cref{it:enum}. Assume that $X$ is r.i.p.e.\ in $\A$. Then there is a computable enumeration of $\Sigma^p_1$ formulas $\psi_n$ with parameters $\bar p$ such that for some $\bar a\in A^{<\omega}$, $n\in X \Leftrightarrow \A\models \psi_n(\bar a)$. Simultaneously enumerate $\Sigma^p_1\text{-}tp_\A(\bar a)$ and the disjuncts in the formulas $\psi_n$. Whenever you see a disjunct of $\psi_n$ that is also in $\Sigma^p_1\text{-}tp_\A(\bar a)$ enumerate $n$. This gives an enumeration of $X$ given an enumeration of $\Sigma^p_1\text{-}tp_\A(\bar a)$ and thus $X\leq_e \Sigma^p_1\text{-}tp_\A(\bar a)$ as required.
\end{proof}

\section{The positive jump and degree spectra}\label{sec:jumps}
In this section we compare various definitions of the jump of a structure with respect to their enumeration degree spectra, a notion first studied by Soskov~\cite{soskov2004}. Before we introduce it, let us recall the definition of the jump of an enumeration degree. 
\begin{defn}
    Let $K_A = \{x \, | \, x \in \Psi^A_x\}$. The \emph{enumeration jump} of a set $A$ is $J_e(A) := A \oplus \overline{K}_A$. The jump of an $e$-degree $\mathbf{a} = \text{deg}_e(A) = \{X : X \equiv_e A\}$ is $\mathbf{a}' = \text{deg}_e(J_e(A))$. We let $\boldsymbol{0}_e = deg_e(\emptyset)$.
\end{defn}
\begin{defn}
    The \emph{positive jump} of a structure $\A$ is the structure 
    \begin{equation*}
        J_p(\A)=(\A, \overline{\vec K^{\A}_p})=(\A,(\overline{K^{\A}_i})_{i\in\omega}).
    \end{equation*}
\end{defn}
We are interested in the degrees of enumerations of $\A$ and $J_p(\A)$. To be precise,
let $f$ be an enumeration of $\omega$, that is, a surjective mapping $\omega \rightarrow \omega$, and for $X\subseteq \omega^{<\omega}$ let
\[ f^{-1}(X)=\{ \langle x_1,\dots\rangle: ( f(x_1),\dots,
)\in X\}.\]
Then, given a structure $\A$ let $f^{-1}(\A)=(\N, f^{-1}(=),f^{-1}(\neq),f^{-1}(R_1^\A),f^{-1}(R_2^\A),\dots)$.
This definition differs from the definition given in Soskov~\cite{soskov2004} where $f^{-1}(\A)$ means what we will denote as $P(f^{-1}(\A))$, i.e.,
\[ f^{-1}(=)\oplus f^{-1}(\neq)\oplus \bigoplus_{i\in\N} f^{-1}(R_i^\A). \]

\begin{defn}[\cite{soskov2004}]\label{def:enumerationspectrum}
  Given a structure $\A$, define the set of enumerations of a structure $Enum(\A)=\{P(\B): \B = f^{-1}(\A) \text{ for } f \text{ an enumeration of } \omega\}$.
  Further, let the \emph{enumeration degree spectrum} of $\A$ be the set
    \[ eSp(\A)=\{ d_e(P(\B)): P(\B)\in Enum(\A)\}\]
  If $\mathbf a$ is the least element of $eSp(\A)$, then $\mathbf a$ is called the
  \emph{enumeration degree} of $\A$.
\end{defn}
In classical computable structure theory many results often do not hold for structures that are structurally too simple, called the automorphically non-trivial structures. One example of this phenomenon is a basic result about Turing degree spectra due to Knight~\cite{knight1986}: The degree spectra of automorphically non-trivial spectra are closed upwards in the Turing degrees. However, if a structure is automorphically trivial, then its degree spectrum contains a single Turing degree. Consider for example the structure $(\mathbb N,=,\neq)$. All its isomorphic copies are computable and in fact, there exists only one copy with universe $\mathbb N$ of this structure.
Considering enumerations of a structure that are not bijective allows us to avoid these dichotomies when studying enumeration degree spectra. It is not hard to see that the enumeration degree spectrum of any structure is closed upwards in the Turing degrees~\cite[Proposition 2.6]{soskov2004}.

As mentioned after \cref{def:cospectra}, Soskov's definition of the co-spectrum of a structure was slightly different~\cite{soskov2004}. He defined the co-spectrum of a structure $\A$ as the set
\[\{\mathbf d: \forall (\mathbf a\in eSp(\A))\, \mathbf d\leq \mathbf a\}. \]
We now show that the two definitions are equivalent.
\begin{proposition}
    For every structure $\A$, $Co(\A)=\{\mathbf d: \forall (\mathbf a\in eSp(\A)) \mathbf d\leq \mathbf a\}$.
\end{proposition}
\begin{proof}
  First note that $\{\mathbf d: \forall (\mathcal C\cong\A) \mathbf d\leq deg_e(P(\mathcal C))\}\supseteq \{\mathbf d: \forall (\mathbf a\in eSp(\A)) \mathbf d\leq \mathbf a\}$ as every $P(\mathcal C)$ is the pullback of $\A$ along an injective enumeration.
  
  On the other hand for any enumeration $f:\omega\to\omega$, $f^{-1}(\A)/f^{-1}(=)\cong\A$. Consider the substructure $\B$ of $f^{-1}(\A)$ consisting of the least element in every $f^{-1}(=)$ equivalence class. Since $f^{-1}(\A)/f^{-1}(=)\cong\A$, we get that $\B\cong \A$. As $P(f^{-1}(\A))$ contains both $f^{-1}(=)$ and $f^{-1}(\neq)$ we can compute an enumeration of $B\oplus \overline B$ from any enumeration of $P(f^{-1}(\A))$ and thus also the graph of its principal function $p_B$. Let $\mathcal C=p_B^{-1}(\B)$. Then $\mathcal C \cong \A$, and $P(\mathcal C)\leq_e P(f^{-1}(\A))$. 
  Thus, 
  \[\{\mathbf d: \forall (\mathbf a\in eSp(\A)) \mathbf d\leq \mathbf a\}=\{\mathbf d: \forall (\mathcal C\cong \A) \mathbf d\leq deg_e(P(\mathcal C))\}=Co(\A).\]
\end{proof}

Using this notation we can see that for any structure $\A$, and any enumeration $f$ of $\omega$ we have that $J_p(f^{-1}(\A))$ is the structure $(f^{-1}(\A),\overline{ \vec K_p^{f^{-1}(\A)}})$. Thus
\begin{equation*}
    P(J_p(f^{-1}(\A))) = f^{-1}(=) \oplus f^{-1}(\neq)\oplus \bigoplus_{j \in \omega}  \overline{K_j^{f^{-1}(\A)}} \oplus\bigoplus_{i \in \omega} R^{f^{-1}(\A)}_i.
\end{equation*}
If we instead apply the enumeration to $J_p(\A)$ we will get the structure $f^{-1}(J_p(\A))$. Now, for every relation $R$ on $\A$, $R^{f^{-1}(\A)} = f^{-1}(R^\A)$ and thus $\overline{ K_j^{f^{-1}(\A)}} = f^{-1}\left(\overline{ K_j^{\A}}\right)$. So $P(f^{-1}(J_p(\A))) = P(J_p(f^{-1}(\A)))$.

\subsection{Properties of R.i.p.e.\ generics}\label{sec:ripegenerics}

Recall \cref{def:ripegeneric} of r.i.p.e.-generic enumerations and presentations. We are now ready to study their properties.
The following lemma is an analogue of the well-known result that there is a $\Delta^0_2$ 1-generic.
\begin{lem}\label{lemma:generic}
    Every structure $\A$ has a r.i.p.e.-generic enumeration $g$ such that $Graph(g) \leq_e P(J_p(\A))$. In particular, $P(g^{-1}(J_p(\A)))\leq_e P(J_p(\A))$.
\end{lem}
\begin{proof}
    Recall our hands-on definition of a r.i.p.e.\ generic enumeration $g$ of $\A$ after \cref{def:ripegeneric} using an enumeration $(R_e)_{e\in\omega}$ of all r.i.p.e.\ subsets of $A^*$.
    The set $\{ \bar p : \exists \bar q \supset \bar p \,, \bar q \in R_e\}$ is $\Sigma_1^p$-definable in $\A$ and so enumerable from $P(J_p(\A))$, which contains $P(\A)$. The coset $\{ \bar p : \forall \bar q \supset \bar p, \bar q \not \in R_e\}$ is co-r.i.p.e. and so enumerable from $\overline{\vec K^\A_p}$. Hence, $P(J_p(\A))$ will be able to decide when a tuple $\bar p_s$ belongs to the upward closure of $R_e$. Thus $Graph(g) \leq_e P(J_p(\A))$.
    
    If we can enumerate the graph of $g$ and also $J_p(\A)$, then to enumerate $g^{-1}(J_p(\A))$, we wait for something of the form $(x,g(x)) \in Graph(g)$ to appear with $g(x) \in J_p(\A)$ and enumerate $x$ into the corresponding slice of $g^{-1}(J_p(\A))$.
\end{proof}
R.i.p.e.\ generic presentations have many useful properties. One of them is that they are minimal in the sense that the only sets enumeration below a r.i.p.e.\ generic presentation are the r.i.p.e. sets.
\begin{lem}\label{generic-ish}
    If $\B$ is a r.i.p.e.-generic presentation of $\A$, then $X\subseteq \omega$ is r.i.p.e. in $\B$ if and only if $X\leq_e P(\B)$.
\end{lem}
\begin{proof}
    If $X$ r.i.p.e. then it is enumerable from $P(\B)$ by definition. Assume $X$ is enumerable from $P(\B)$. Then $X = \Psi^{P(\B)}_e$ for some $e$. Recall the set $Q_e$ from \cref{ripe} which we know to be r.i.p.e. because we gave a $\Sigma_1^p$ description of it.
    \begin{equation*}
    Q_e = \{ \overline{q} \in A^{< \omega} : \exists l, j_1 ,\dots j_l < |\overline{q}| \, [\langle j_1,\dots, j_l \rangle \in \Psi_e^{P_\B(\overline{q})} \text{ and } \langle q_{j_1},\dots, q_{j_l} \rangle \not \in X] \}.
    \end{equation*}
    Now because $\B$ is r.i.p.e.-generic, there is some tuple $\langle 0,\dots ,k-1 \rangle$ that decides $Q_e$. It must be that $\langle 0,\dots ,k-1 \rangle \not \in Q_e$, or we would contradict the fact that $X = \Psi^{P(\B)}_e$, and so $\langle 0,\dots ,k-1 \rangle$ forms the parameterizing tuple $\bar p_s$ from the set $S$ in \cref{ripe}.
\end{proof}

Another useful property is that the degree of the positive jump of a r.i.p.e.\ generic agrees with the enumeration jump of the degree of their positive diagram.
\begin{proposition}\label{propn generic}
Let $\A$ be a structure. For an arbitrary enumeration $f$ of $\A$, let $\B = f^{-1}(\A)$. Then $P(J_p(\B)) \leq_e J_e(P(\B))$. Furthermore, if $\B$ is a r.i.p.e.-generic presentation then $P(J_p(\B)) \equiv_e J_e(P(\B))$.
\end{proposition}
\begin{proof}
Recall that $J_p(\B) = (\B, \overline{\vec K_p^{\B}})$ and $J_e(P(\B)) = P(\B) \oplus \overline{K}_{P(\B)}$, so to show that $P(J_p(\B)) \leq_e J_e(P(\B))$ it is sufficient to show that $\overline{\vec K_p^{\B}} \leq_e \overline{K}_{P(\B)}$.
As $\vec K_p^{\B}$ is r.i.p.e., $\vec K_p^{\B} = \Psi^{P(\B)}_e$ for some $e$. Thus $\overline{\vec K_p^{\B}}$ appears as a slice of $\{\langle x,i \rangle : x \not \in \Psi^{P(\B)}_i\}$ and hence $\overline{\vec K_p^{\B}}\leq_e\{\langle x,i \rangle : x \not \in \Psi^{P(\B)}_i\}$. The latter set is $m$-equivalent to $\overline{K}_{P(\B)}$ and thus $\overline{\vec K_p^{\B}}\leq_e\overline{K}_{P(\B)}$.

It remains to show that $J_e(P(\B))\leq_e P(J_p(\B))$ if $\B$ is r.i.p.e.-generic. For every $e$ consider the set
\[ R_e = \{\bar q\in B^* : e \in \Psi_e^{P_\B(\bar q)}\}.\]
Note that the $R_e$ are all trivially r.i.p.e. and are closed upwards as subsets of $B^*$. So since $\B$ is r.i.p.e.-generic, for every $e$ there is an initial segment of $B^*$ that either is in $R_e$ or such that no extension of it is in $R_e$. The set $Q_e=\{ \bar p\in B^* :\forall (\bar q\supseteq \bar p) (\bar q \not \in R_e)\}$ of elements in $B^*$ that are non-extendible in $R_e$ is co-r.i.p.e. Also note that these sets are uniform in $e$, that is, given $e$ we can compute the indices of $R_e$ and $Q_e$ as r.i.p.e., respectively, co-r.i.p.e.\ subsets of $\B$. 
Thus, given an enumeration of $P(J_p(\B))$ we can enumerate $P(\B)$ and the sets $Q_e$ and $R_e$. By genericity for all $e\in \omega$ there is an initial segment of $\bar b$ of $B$ such that either $\bar b\in Q_e$ or $\bar b \in R_e$. So whenever we see such a $\bar b$ in $Q_e$ we enumerate $e$ and thus obtain an enumeration $\overline K_{P(\B)}$.
\end{proof} 
The above properties of r.i.p.e.\ generics are very useful to study how the enumeration degree spectra of structures and their positive jumps relate.

\begin{proposition}
For every structure $\A$, $J_p(\A)$ is a total structure.
\end{proposition}
\begin{proof}
    First, note that we have the following equality for arbitrary relations $R$ and enumerations $f$.
    \[ x\in f^{-1}(\bar R) \Leftrightarrow f(x)\in \bar R \Leftrightarrow f(x)\not\in R \Leftrightarrow x\not\in f^{-1}(R)\Leftrightarrow x\in \bar{f^{-1}(R)}\]
    Recall that for every $R_i$ in the language of $\A$, $R_i$ is r.i.p.e. in $\A$ and that 
    \[ P(f^{-1}(J_p(\A))) = f^{-1}(=) \oplus f^{-1}(\neq) \oplus f^{-1}(\overline{\vec K^{\A}_p}) \oplus f^{-1}(R^\A_1) \oplus \cdots\]
    and we observe that all $R_i^\A$ are trivially r.i.p.e., uniformly in $i$, so in particular $\overline{f^{-1}(R_i^\A)} \leq_e \overline{f^{-1}(\vec K^\A_p)}=f^{-1}(\bar{\vec K^\A_p})$. Also, $f^{-1}(\vec K^\A_p)\leq_e P(f^{-1}(\A))\leq_e P(f^{-1}(J_p(\A)))$.
\end{proof}
We will now see that the positive jump of a structure jumps, in the sense that the enumeration degree spectrum of the positive jump is indeed the set of jumps of the degrees in the enumeration degree spectrum of the structure. The following version of a Theorem by Soskova and Soskov~\cite{soskova2009} is essential to our proof.
\begin{thm}[{\cite[Theorem 1.2]{soskova2009}}]\label{Soskov}
    Let $B$ be an arbitrary set of natural numbers. There exists a total set $F$ such that
    \begin{equation*}
        B \leq_e F \quad \text{and} \quad J_e(B) \equiv_e J_e(F).
    \end{equation*}
\end{thm}

\begin{thm}\label{theorem_espjumps}
    For any structure $\A$,
    \begin{equation*}
        eSp(J_p(\A)) = \{ \mathbf{a'} : \mathbf{a} \in eSp(\A) \}
    \end{equation*}
\end{thm}

\begin{proof}
    To show that $\{ \mathbf{a'} : \mathbf{a} \in eSp(\A) \}\subseteq
    eSp(J_p(\A))$, consider an arbitrary enumeration $f$ of $\omega$ and let $\B
    = f^{-1}(\A)$. Then from \cref{propn generic} we know that $P(J_p(\B)) \leq_e J_e(P(\B))$. Note that $P(J_p(\B))=P(J_p(f^{-1}(\A))) = P(f^{-1}(J_p(\A)))$, and $d_e(P(f^{-1}(J_p(\A))) \in eSp(J_p(\A))$. Since the enumeration jump is always total and enumeration degree spectra are closed upwards with respect to total degrees, $d_e(J_e(P(\B))) \in eSp(J_p(\B))$. 
    
    To show $eSp(J_p(\A)) \subseteq \{ \mathbf{a'} : \mathbf{a} \in eSp(\A) \}$, let us look at $P(f^{-1}(J_p(\A)))$ for some enumeration $f$ of $\omega$, and again write $\B = f^{-1}(\A)$ so that $P(f^{-1}(J_p(\A))) = P(J_p(\B))$. Since $J_p(\A)$ is a total structure, we know that $P(J_p(\B))$ is total. By \cref{lemma:generic} we can use $P(J_p(\B))$ to enumerate a r.i.p.e.-generic enumeration $g$ of $\B$ such that $P(g^{-1}(J_p(\B))) \leq_e P(J_p(\B))$. Then, letting $\mathcal{C} = g^{-1}(\B)$ and using the latter part of \cref{propn generic} we know that $J_e(P(\mathcal{C})) \equiv_e P(J_p(\mathcal{C}))$ which makes $J_e(P(\mathcal{C})) \leq_e P(J_p(\B))$. Using \cref{Soskov} there is a total set $F$ such that $P(\mathcal{C}) \leq_e F$ and $J_e(F) \equiv_e J_e(P(\mathcal{C})) \leq_e P(J_p(\B))$.
    As $F$ and $P(J_p(\B))$ are total, we have that $F' \leq_T P(J_p(\B))$.
    
    We can now apply the relativized jump inversion theorem for the Turing degrees to get a set $Z\geq_T F$ such that $Z' \equiv_T P(J_p(\B))$. For this $Z$, we have that $\iota(deg(Z')) = d_e(J_e(\chi_Z)) = d_e(P(J_p((\B))) \in eSp(J_p(\A))$ and $d_e(P(\B)) \leq d_e(F) \leq d_e(\chi_Z)$. Since $\chi_Z$ is total and enumeration degree spectra are upwards closed in the total degrees, this means $d_e(\chi_Z) \in eSp(\B) \subseteq eSp(\A)$. So in particular, $d_e(J_e(\chi_Z))=d_e(P(J_p(\B)))\in \{\mathbf a': \mathbf a\in eSp(\A)\}$.
\end{proof}
We now compare the enumeration degree spectrum of the positive jump of a structure $\A$ with the spectrum of the original structure and the spectrum of its traditional jump. The latter is the spectrum of the structure obtained by adding the sequence of relations
\[  \vec K^{\A} = (\bigcup_{j\in\omega}\{ \bar a\in \omega^j : \A\models \phi_{i,j}^{\Sigma_1^c}(\bar a)\})_{i\in\omega} \]
where $\phi_{i,j}^{\Sigma_1^c}(\bar x)$ is the $i^{th}$ formula with $j$ free variables in a fixed computable enumeration of all computable infinitary $\Sigma_1^0$ formulas.
Let $J(\A)=(\A,\vec K^{\A})$ denote the traditional jump of $\A$. This notion of a jump was defined by Montalb\'an~\cite{montalban2012} and is similar to several other notions of jumps of structures that arose independently such as Soskova's notion using Marker extensions~\cite{soskova2009} or Stukachev's version using $\Sigma$-definability~\cite{stukachev2007}.  

We will also consider the structure $\A^+$, the totalization of $\A$ given by
\[ \A^+=(A,(R_i^\A,\bar{ R_i^\A})_{i\in\omega}).\]

We will not compare the enumeration degree spectra directly but instead the sets of enumerations. This gives more insight than comparing the degree spectra as we can make use of the following analogues to Muchnik and Medvedev reducibility for enumeration degrees. Given $\mf A,\mf B\subseteq P(\N)$ we say that $\mf A\leq_{we} \mf B$, $\mf A$ is \emph{weakly enumeration reducible} to $\mf B$, if for every $B\in \mf B$ there is $A\in \mf A$ such that $A\leq_e B$. We say that $\mf A\leq_{se}\mf B$, $\mf A$ is \emph{strongly enumeration reducible} to $\mf B$, if there is an enumeration operator $\Psi$ such that for every $B\in\mf B$, $\Psi^B\in \mf A$.

It is not hard to see that given an enumeration of $\A^+$ one can enumerate an enumeration of $J(\A)$ and the converse holds trivially. We thus have the following.
\begin{proposition}\label{prop:tradjumpiscompletion}
    For every structure $\A$, $deg_e(J(\A))=deg_e(\A^+)$. In particular $Enum(J(\A))\equiv_{se} Enum(\A^+)$.
\end{proposition}
\begin{proof}
By replacing every subformula of the form $\neg R_i(x_1,\dots, x_m)$ by the formula $\bar R_i(x_1,\dots,x_m)$ we get a $\Sigma^p_1$ formula in the language of $\A^+$. Similarly, given a $\Sigma^p_1$ formula in the language of $\A^+$ we can obtain a $\Sigma^c_1$ formula in the language of $\A$ by substituting subformulas of the form $\bar R_i(x_1,\dots,x_m)$ with $\neg R_i(x_1,\dots,x_m)$. Indeed we get a computable bijection between the $\Sigma^p_1$ formulas in the language of $\A^+$ and the $\Sigma^c_1$ formulas of $\A$. Thus $J(\A)\equiv_e (\A^+,\vec K_p^{\A^+})$, but $\A^+\equiv_e (\A^+,\vec K_p^{\A^+})$ given the first equivalence. All of these equivalence are witnessed by fixed enumeration operators and thus $Enum(J(\A))\equiv_{se} Enum(\A^+)$.
\end{proof}
\cref{prop:tradjumpiscompletion} is not very surprising, as adding a r.i.c.e.\ set such as $\vec{K}^\A$ to the totalization of $\A$ won't change its enumeration degree. We will thus consider a slightly different notion of jump by adding the coset of $\vec{K}^\A$:
\[ T(\A)=(\A,\overline{\vec{K}^{\A}})\]
Clearly $T(\A)\equiv_T J(\A)$; however, the two sets are not necessarily enumeration equivalent as the following two propositions show.
\begin{proposition}
    Let $\A$ be a structure. For every presentation $\hat \A$ of $\A$, $P(\hat \A)\leq_e D(\hat\A)=P(\hat\A^+)\leq_e J_p(\hat\A)\leq_e T(\hat\A)$. In particular 
    \[Enum(\A)\leq_{se}Enum(\A^+)\leq_{se} Enum(J_p(\A))\leq_{se} Enum(T(\A)).\]
\end{proposition}
\begin{proof}
Straightforward from the definitions.
\end{proof}
\begin{proposition}
There is a structure $\A$ such that \[Enum(\A)\not\geq_{we}Enum(\A^+)\not\geq_{we} Enum(J_p(\A))\not\geq_{we}Enum(T(\A)).\]
\end{proposition}
\begin{proof}
    Sacks showed that there is an incomplete c.e.\ set $X$ of high Turing degree. Thus, $X$ has enumeration degree $deg_e(X)=\mathbf 0_e$ and $deg_T(X')=\mathbf 0''$. Let $\A$ be the graph coding $X$ as follows. $\A$ contains a single element $a$ with a loop, i.e. $aE^{\A}a$ and one circle of length $n+1$ for every $n\in\omega$. Let $y$ be the least element in $\A$ that is part of the cycle of length $n+1$. If $n\in X$, then connect $a$ to $y$, i.e., $aE^\A y$. This finishes the construction. 
    
    As $X$ is c.e., $\A$ has enumeration degree $\mathbf 0_e$. However, $\mbf 0_e\not\in eSp(\A^+)$, as $\A^+$ has enumeration degree $deg_e(X\oplus\overline X)> \mathbf 0_e$. So, in particular $Enum(\A)\not\geq_{we} Enum(\A^+)$.
    By \cref{theorem_espjumps} the enumeration degree of $J_p(\A)$ is $\mathbf 0'_e$ and $\mathbf 0'_e\ni \emptyset'\oplus \overline{\emptyset'}>_e X\oplus\overline X$, so $Enum(\A^+)\not\geq_{we}Enum(J_p(\A)) $. 
    For the last inequality notice that both $T(\A)$ and $J_p(\A)$ are total structures and that by the analogue of \cref{theorem_espjumps} for the traditional jumps of structures we have that the enumeration degree of $T(\A)$ is  $deg_e(\emptyset''\oplus \overline{\emptyset''})$. As mentioned above the enumeration degree of $J_p(\A)$ is $\mathbf 0'_e$, $\emptyset'\oplus \overline{\emptyset'}\in \mathbf 0'_e$, and $\emptyset'\oplus \overline{\emptyset'}<_e \emptyset''\oplus \overline{\emptyset''}$. So, in particular, $Enum(J_p(\A))\not \geq_{we} Enum(T(\A))$.
\end{proof}
For total structures the traditional notion of the jump and the positive jump coincide.
\begin{proposition}
    Let $\A$ be a structure, then $Enum(J_p(\A^+))\equiv_{se}Enum(T(\A))$ and $eSp(J_p(\A^+))=eSp(T(\A))$.
\end{proposition}
\begin{proof}
    This proof is similar to the proof of \cref{prop:tradjumpiscompletion}. Mutatis mutandis.
\end{proof}

\section{Conclusion and open questions}
In an earlier version of this paper, we claimed a theorem showing that a functorial notion of reduction introduced in~\cite{CsimaRY21} is equivalent to a syntactic notion of interpretability using $\Sigma^p_1$ formulas. Unfortunately, our proof was wrong and all efforts to fix it failed. In the Turing computability setting, Harrison-Trainor, Melnikov, Miller, and Montalbán~\cite{HTMMM2015} showed that computable functors are equivalent to effective interpretability---a version of classical model-theoretic interpretability that allows definitions via $\Sigma^c_1$ formulas. As the notion of a positive enumerable functor introduced in~\cite{CsimaRY21} preserves the enumeration analogues of properties preserved by computable functors, we still believe that there should be a version of interpretability using $\Sigma^p_1$ formulas that yields a syntactic equivalent. We thus pose the following question.
\begin{question}
Is there a notion of interpretability, in particular, one using $\Sigma^p_1$ formulas, such that a structure $\A$ is interpretable in $\B$ if and only if there is a positive enumerable functor $F:Iso(\B)\to Iso(\A)$?
\end{question}
Since the first version of this article, other work relating enumeration degrees and computable structure theory has emerged. In~\cite{bazhenov2023} a Lopez-Escobar theorem for Scott topologies was proven: It was shown that an isomorphism invariant class of structures is $\Sigma^p_\alpha$-definable if and only if the class is $\pmb\Sigma^0_\alpha$ in the Borel hierarchy of the Scott topology on this space. This result suggests that it should be possible to generalize some of our results on r.i.p.e.\ relations through the transfinite.
\printbibliography

\end{document}